\documentclass[11pt]{article}
\usepackage{amsmath,amsthm,amsfonts,latexsym,amssymb,enumerate,color}

\def\CC{\mathbb{C}}
\def\DD{\mathbb{D}}

\def\GG{\mathcal{G}}
\def\KK{\mathcal{K}}

\def\NN{\mathbb{N}}
\def\RR{\mathbb{R}}

\def\ZZ{\mathbb{Z}}

\def\Hp{H_p^+}

\def\Np{\widetilde{\mathcal{N}_p}}

\def\dim{\mathop{\rm dim}\nolimits}
\def\ker{\mathop{\rm ker}\nolimits}
\def\Kmin#1{\KK_{\rm min}\left( #1 \right) }
\def\LCM{\mathop{\rm LCM}\nolimits}

\def\nei#1{n.~$#1$-in\-vari\-ant}

\renewcommand\span{\mathop{\rm span}\nolimits}
\def\spam{\span} % I always mistype this
\def\clos{{\rm clos}}

\def\ds{\displaystyle}
\def\ol{\overline}

\def\th{\theta}

\def\eps{\varepsilon}

\renewcommand\phi{\varphi}

\newtheorem{thm}{Theorem}[section]
\newtheorem{prop}[thm]{Proposition}
\newtheorem{defn}[thm]{Definition}

\newtheorem{cor}[thm]{Corollary}
\newtheorem{rem}[thm]{Remark}

\newtheorem{ex}[thm]{Example}

\numberwithin{equation}{section}

\def\beginpf{\begin{proof}}
\def\endpf{\end{proof}}
\def\beq{\begin{equation}}
\def\eeq{\end{equation}}

\def\imag{\mathop{\rm Im}\nolimits}

\begin{document}

\title{Model spaces and Toeplitz kernels in reflexive Hardy spaces}
\author{M.~C. C\^amara\thanks{
Center for Mathematical Analysis, Geometry and Dynamical Systems,
Mathematics Department,
Instituto Superior T\'ecnico, Universidade de Lisboa,
Av. Rovisco Pais, 1049-001 Lisboa, Portugal.
 \tt ccamara@math.ist.utl.pt}~, 
M.T. Malheiro\thanks{
Centre of Mathematics, Departamento de Matem\'atica e Aplica\c{c}\~oes, Universidade do Minho, Campus de Azur\'em, 4800-058 Guimar\~aes, Portugal. \tt mtm@math.uminho.pt}
~and J. R.~Partington\thanks{School of Mathematics, University of Leeds, Leeds LS2~9JT, U.K. \tt j.r.partington@leeds.ac.uk}
}

\maketitle

\begin{abstract}
This paper considers model spaces  in an $H_p$ setting. 
The existence of unbounded functions and the characterisation of maximal functions in a model space are studied, and decomposition results for Toeplitz kernels, in terms of model spaces, are established.
\end{abstract}

%%%%%%%%%%%%%%%%%%%%%%%%%%%%%%%                            Section 1            Introduction                                               %%%%%%%%%%%%%%%%%%%%%%%%%%%%%%%%%%%%%%%%%%%
%%%%%%%%%%%%%%%%%%%%%%%%%%%%%%%%%%%%%%%%%%%%%%%%%%%%%%%%%%%%%%%%%%%%%%%%%%%%%%%%%%%%%%%%%%%%%%%%%%%%%%%%%%%%%%%%%%

\section{Introduction and notation}

In the theory of complex functions and linear operators, there has been a significant body of
work attempting to understand the structure and properties of kernels of Toeplitz
operators, or Toeplitz kernels, and to describe them (or at least determine their dimension) explicitly for some concrete classes of symbols (see, for example, \cite{BCD,CP14,Dyakonov,hayashi86,sarason88,sarason}).\\

Linked with this is the theory of {\em model spaces\/}, which have generated an enormous interest; they provide the natural setting for truncated Toeplitz operators and are relevant in connection with the study of a variety of topics such as the Schr\"odinger operator, classical extremal problems, and Hankel operators 
(see for instance \cite{Garcia and Ross} and references therein).\\

Model spaces constitute a particular type of Toeplitz kernel whose properties are in general more fully understood.
Indeed, denoting by $\mathbb{D}$ the unit disk, Beurling's theorem characterises the nontrivial subspaces of $H^2(\mathbb D)$ which are invariant under the (unilateral) shift $S$ as consisting of the $H^2(\mathbb D)$ multiples of some inner function $\theta$, i.e., as being of the form $\theta H^2(\mathbb D)$. The so-called model spaces $K_\theta$ are the nontrivial invariant subspaces for the backward shift $S^*$; they are the orthogonal complements in $H^2(\mathbb D)$ of the shift-invariant subspaces $\theta H^2(\mathbb D)$.\\

An equivalent definition, which is better suited to   the context of the Hardy spaces $H^p$ with $p\in (1,\infty), p\neq 2$, in which the Hilbert space structure is absent, is to say that $K_\theta$ is the kernel of the Toeplitz operator whose symbol is $\ol \theta$, the complex conjugate of the inner function $\theta$, assumed to be non-constant.
This approach to model spaces in $H^p(\mathbb D)$, or in $H_p^+:= H^p(\mathbb C^+)$ which will be our main setting (here $\CC^+$ denotes the upper half-plane), provides a simple operator theory point of view, as well as a functional analytic description of $S^*$-invariant subspaces which is almost as simple as Beurling's description of $S$-invariant subspaces: $K_\theta$ consists of the $H_p^-$ multiples of $\theta$ which belong to $H_p^+$ (using the notation $H_p^\pm$ for $H^p({\mathbb C}^\pm))$.\\

This paper's results
take further some ideas introduced
in \cite{CP14}, looking at model spaces and Toeplitz
operators in a general $H_p$ context ($1<p<\infty$),
rather than simply $H_2$, and working on the upper half-plane rather than the disk. One advantage of this choice is that some formulae are simpler in the half-plane context, although they can generally be translated to
analogous results on the disk; some questions, however, are meaningful only in a half-plane context.\\

The themes considered in this work include near invariance (a property of all Toeplitz kernels, and model spaces in particular), the dependence of a Toeplitz kernel on the symbol of the corresponding Toeplitz operator and the $H_p$ space where it is defined, some associated factorisation and decomposition results, and the
existence of a {\em maximal function\/} in every Toeplitz kernel that uniquely defines the latter.
The results also generalise some properties of model spaces to general Toeplitz kernels and show that we can use model spaces to ``quantify" (in a loose sense of the word), for infinite-dimensional kernels, some properties relating the dimensions of finite-dimensional kernels.\\

More precisely, the structure of this paper is as follows. The first two sections are of an auxiliary nature. In Section~\ref{sec:2} we present some results
on Toeplitz kernels and near invariance in an $H_p^+$ context, very much in the spirit of \cite{CP14}.
In Section~\ref {sec:3} we turn our attention to model spaces, regarded as
Toeplitz kernels of a particular kind, and present their basic properties and some factorisation and decomposition
results. The main results of the paper are contained in the next three sections. Section~\ref{sec:4} addresses the question when model spaces consist entirely of
bounded functions, i.e., form subspaces of $H_\infty^+$; the answer for the
half-plane turns out to be significantly more interesting than in the disk case and provides an example where results on the disk do not carry over to the upper half-plane and vice-versa.
Then in Section~\ref{sec:5} we are mainly concerned with characterising maximal functions in a model space, i.e., those which are
contained in no smaller Toeplitz kernel.
Finally, in Section~\ref{sec:6} we  establish decomposition results relating two Toeplitz kernels 
determined by symbols that differ only by an inner factor.\\

We take $1<p<\infty$ and $H_p^+$, $H_p^-$ to be the Hardy spaces of the upper and lower half-planes $\CC^+$ and $\CC^-$ respectively.
We write $L_p$ to denote $L^p(\RR)$. 

The class of invertible elements in $H_\infty^\pm$ is denoted by $\GG H_\infty^\pm$. Similarly for $\GG L_\infty$.

We write $P^+:L_p  \to H^+_p$ for the projection
with kernel $H^-_p$.\\

For $g \in L_\infty(\RR)$  and $1<p<\infty$,  the {\em Toeplitz operator\/}
$T_g: H^+_p \to H^+_p$ is defined by
\[
T_g f_+ = P^+(g f_+), \qquad (f_+ \in H^+_p).
\]
We shall require the functions
\beq \label {eq:1.1}
\lambda_{\pm}(\xi)= \xi \pm i \qquad \hbox{and} \qquad r(\xi)= \frac{\xi-i}{\xi+i}\,\,,
\eeq
 and write $S$ for the operator $T_r$ on $\Hp$  of multiplication by $r$, with $S^*$ the
operator $T_{\overline r}$.

%-------------------------------------------------------------------------------------------------SECTION 2-------

\section{Near invariance and T-kernels}
\label{sec:2}

%\begin{defn}
%Let $\E$ be a proper closed subspace of {$\Hp$}, and $\eta $ a complex-valued function defined a.e.\ on $\RR$. We say that $\E$ is
%{\em nearly $\eta$-invariant}, if for every $f_+ \in \E$ such that $\eta f_+ \in \Hp$, we have $\eta f_+ \in \E$; that is,
%$\eta\E \cap \Hp \subset \E$.
%If $\E$ is nearly $\eta$-invariant with $\eta\in L_\infty$, then we also say that $\E$ is {\em nearly $T_\eta$-invariant.}
%\end{defn}
%
%We abbreviate ``nearly $\eta$-invariant'' to ``\nei{\eta}".\\
%
%Defining $\mathcal{N}_p$ as the set of all complex-valued functions $\eta$, defined a.e.\ on $\mathbb{R}$, such that every kernel of a Toeplitz operator (abbreviated to T-kernel) in $H_p^+$ is \nei{\eta},
%it is shown in \cite{CP14}  that
%$$\Np\subset\mathcal{N}_p,$$
%where  $\Np :=\{\eta: L_p \cap \eta H_p^- \subset H_p^- \}$.\\

%Definition 2.1

\begin{defn}  \cite{CP14}
 Let  $\mathcal{E}$ be a proper closed subspace of $H_p^+$ and $\eta$ a complex-valued function defined almost everywhere on $\RR$. We say that $\mathcal{E}$ is \emph{nearly}~$\eta$-\emph{invariant} if and only if, for every $f_+\in \mathcal{E}$ such that $\eta f_+ \in H_p^+$, we have $\eta f_+ \in \mathcal{E}$; that is
\begin{equation}\label{21.1}
\eta\mathcal{E} \cap H_p^+ \subset \mathcal{E} .
\end{equation}
If $\mathcal{E}$ is nearly $\eta$ - invariant with $\eta \in L_{\infty}$, then we also say that $\mathcal{E}$ is \emph{nearly}~$T_\eta$-\emph{invariant}.
\end{defn}

We abbreviate ``nearly $\eta$ - invariant" to ``\nei{\eta}".\\ 

We denote by $\mathcal{N}_p$ the set of all complex-valued functions $\eta$, defined a.e. on $\RR$, such that every kernel of a Toeplitz operator (abbreviated to \emph {T-kernel}) in $H_p^+$ is \nei{\eta}, i.e., such that for all $g \in L_\infty$ we have
\begin{equation}\label{21.2}
\eta\ker T_g\cap H_p^+ \subset \ker T_g.
\end{equation}

It is shown in \cite{CP14} that $\mathcal{N}_p \supset \widetilde{\mathcal{N}}_p$, where
\[\widetilde{\mathcal{N}}_p :=\{\eta:L_p \cap \eta H_p^-\subset H_p^-\},\]
and that many well-known classes of functions are contained in $\widetilde{\mathcal{N}}_p$, amongst them 
$\mathcal{L}^-_{\infty,m}:=\lambda_-^m H_\infty^-$ for all $m \in \mathbb{Z}$, the set of all rational functions with poles 
belonging to $\mathbb{C}^+\cup \RR \cup \{\infty\}$, and $H_p^-$ for all $p \in \,(1, \infty)$.\\

On the other hand, if we extend the notation for T-kernels, defining
\begin{equation}\label{21.3}
\ker T_g:=\{\phi_+ \in H_p^+: g\phi_+ \in H_p^-\}
\end{equation}
for all complex-valued $g$ defined a.e. on $\RR$, it is clear that we also have
\begin{equation}\label{21.4}
\eta\ker T_g\cap H_p^+ \subset \ker T_{\eta^{-1}g}
\end{equation}
if $\eta^{\pm 1}$ are defined a.e. on $\RR$ (whether or not they belong to $\mathcal{N}_p$). We have moreover:

%Proposition 2.2
\begin{prop}\label{thm:2.2} If $\eta \in \widetilde{\mathcal{N}}_p$, then $\ker T_{\eta^{-1}g}\subset \ker T_g$ for all $g \in L_\infty$.
\end{prop}

\beginpf Let $\phi_+ \in H_p^+$ and $\eta^{-1}g\phi_+=\phi_- \in H_p^-$. Then $g\phi_+=\eta\phi_- \in L_p \cap \eta H_p^- \subset H_p^-$, so that $\phi_+ \in \ker T_g$.
\endpf

Taking \eqref{21.4} into account we have thus:

%Corollary 2.3
\begin{cor}\label{cor:2.3} If $\eta \in \widetilde{\mathcal{N}}_p$, $g \in L_\infty$, then
\begin{equation}\label{21.5}
\eta \ker T_g \cap H_p^+ \subset \ker T_{\eta^{-1}g} \subset \ker T_g.
\end{equation}
\end{cor}

The inclusions in \eqref{21.4} and in Proposition \ref{thm:2.2} may be strict or not.

Regarding the first inclusion, it is easy to see that if $ {O}_+$ is outer in $H_\infty^+ $ then

\begin{equation}\label{21.6}
{O}_+^{-1} \ker T_g \cap H_p^+ = \ker T_{ {O}_+ g},
\end{equation}

and in particular

\begin{equation}\label{21.7}
h_+\in \mathcal{G}H_\infty^+ \Rightarrow h_+\ker T_g = \ker T_{h_+^{-1} g}.
\end{equation}

On the other hand, for any non-constant inner function $\theta$, if $\ker T_g \neq \{0\}$ then

\begin{equation}\label{21.8}
\theta\ker T_{\theta g} \varsubsetneq \ker T_{g}
\end{equation}

since either $\ker T_{\theta g}= \{0\}$ and  \eqref {21.8} is obvious, or $\ker T_{\theta g}\neq \{0\}$ and \eqref {21.8} 
follows from \eqref {21.4} and the proposition below.

%Proposition 2.4
\begin{prop}\label{thm:2.4} If $\tilde g$ is a complex-valued function defined a.e. on $\RR$, $\ker T_{\tilde g}\neq \{0\}$ and $\theta$ is a non-constant inner function, then $\theta \ker T_{\tilde g}$ is not a \nei{\bar{\theta}} subspace of $H_p^+$.
\end{prop}

\beginpf For $\mathcal{E}= \theta \ker T_{\tilde g}$, we have $\bar{\theta} \mathcal{E}= \ker T_{\tilde g} \subset H_p^+$. But if $\ker T_{\tilde g} \subset \mathcal{E}$, then for any $\phi_+ \in \ker T_{\tilde g}$ we would have $\phi_+= \theta \psi_+$ with $\psi_+ \in \ker T_{\tilde g}$ and, repeating this reasoning, $\phi_+$ would be divisible in $H_p^+$ by arbitrarily large powers of $\theta$, implying that $\phi_+=0$.
\endpf

We remark however that $\theta \ker T_{\tilde g}$ is a \nei{S^*} subspace of $H_p^+$ if $\theta (i)\neq 0$. Indeed if $r^{-1}\theta \phi_+ \in H_p^+$, with $\phi_+ \in \ker T_{\tilde g}$, then we must have $\phi_+(i)=0$, so that $r^{-1}\phi_+ \in H_p^+$, and ${\tilde g}r^{-1}\phi_+=r^{-1}\phi_-$ with $\phi_- \in H_p^-$, implying that $r^{-1} \phi_+ \in \ker T_{\tilde g}$ and $r^{-1} \theta \phi_+\in \theta \ker T_{\tilde g}$.\\

Regarding the inclusion in Proposition \ref{21.2}, we have the following two results.

% Prop 2.5
\begin{prop} If $\eta^{\pm 1} \in \widetilde{\mathcal{N}}_p$ and $g,\, \eta g \in L_\infty$, then $\ker T_{\eta^{-1}g}=\ker T_g$.
\end{prop}

\beginpf From Corollary \ref{cor:2.3} we have, on the one hand, $\ker T_{\eta^{-1}g}\subset\ker T_g$ and, on the other hand, $\ker T_g=\ker T_{\eta (\eta^{-1}g)}\subset \ker T_{\eta^{-1}g}$.
\endpf

In particular, if $ {O}_{-}$ is outer in $H_\infty^-$ then
\begin{equation}\label{21.9}
\ker T_{ {O}_-g} = \ker T_g
\end{equation}
and
\begin{equation}\label{21.10}
h_- \in \mathcal{G}H_\infty^- \Rightarrow \ker T_{h_-g}=\ker T_g.
\end{equation}

%Proposition 2.6
\begin{prop} If $\eta= \bar{\theta} \tilde{\eta}$, where $\theta\in H_\infty^+$ is a non-constant inner function and $\tilde{\eta} \in \widetilde{\mathcal{N}}_p$, then
\[\ker T_{\eta^{-1}g}\varsubsetneq \ker T_g\]
if $\ker T_g \neq \{0\}$.
\end{prop}

\beginpf If $\ker T_{\eta^{-1}g}=\{0\}$, the inclusion is obviously strict. If $\ker T_{\eta^{-1}g}\neq\{0\}$ then, by an analogue of Theorem 2.2 in \cite{CP14} and Proposition \ref{thm:2.2} above,
\[\ker T_{\eta^{-1}g}=\ker T_{\theta\tilde{\eta}^{-1}g}\varsubsetneq \ker T_{\tilde{\eta}^{-1}g} \subset \ker T_g.\]
\endpf

Note that studying T-kernels is closely related to studying sets of the form $\eta \ker T_g \cap H_p^+$ since we can write, for the kernel of any operator $T_G$ in $H_p^+$,
\begin{equation}\label{21.11}
\ker T_G=h_+ (\bar{\theta}_1\ker T_{\bar{\theta}_2} \cap H_p^+)
\end{equation}
where $h_+ \in \mathcal{G}H_\infty^+$ and $\theta_1$, $\theta_2$ are inner functions, which may be chosen to be Blaschke products (\cite{Dyakonov}, Theorem 1).

%%%%%%%%%%%%%%%%%%%%%%%%%%%%%%%%%%        Section 3       Model spaces in Hp         %%%%%%%%%%%%%%%%%%%%%%%%%%%%%%%%%%%%%%%%%%%
%%%%%%%%%%%%%%%%%%%%%%%%%%%%%%%%%%%%%%%%%%%%%%%%%%%%%%%%%%%%%%%%%%%%%%%%%%%%%%%%%%%%%%%%%%%%%%%%%%%%%

\section{Model spaces in $H_p^+$}\label{sec:modelspaces}
\label{sec:3}

%Definition 3.1
\begin{defn}
If $\theta$ is an inner function, then $K_{\theta}^p:= H_p^+ \cap \, \theta H_p^-$, for $p \in (1, \infty)$.
\end{defn}

We omit the superscript $p$ in $K_{\theta}^p$ unless it is required for clarity.\\

This definition makes it clear that $K_\theta$ is a T-kernel, since $K_\theta= \ker T_{\ol\theta}$.
Model spaces are thus \nei{\eta} for all $\eta \in \Np$; in the case of $\eta \in H_\infty^-$, model spaces are moreover $T_{\eta}$-invariant.
A particular case is that of $S^*=T_{r^{-1}}$, where $r$ is given by \eqref {eq:1.1},
in which case the converse is true (\cite{Gamelin})  and we can say that $K \subset H_p^+$ is a model space
if and only if $K$ is $S^*-invariant$.\\

%{\bf This is easily deduced from known results, but I need to give proper explanations. JRP}\\

Given $p\in (1,\infty)$, to each inner function $\theta$ we can associate a bounded projection $P_\theta: L_p \to K_{\theta}$ defined by
\begin{equation}\label{2.2}
P_\theta = \theta P^- \overline\theta P^+.
\end{equation}
Its restriction to $H_p^+$ is also a projection onto $K_{\theta}$, which we denote in the same way.
We have $K_{\theta}=P_\theta H_p^+ = P_\theta L_p$ and $H_p^+ = K_{\theta} \oplus \theta H_p^+$ (for $p=2$ this is
an orthogonal decomposition).\\

We also have $K_{\theta}=P^+(\theta H_p^-)$ and
\begin{equation}\label{2.3}
K_\theta=\theta \, \overline{K_\theta}.
\end{equation}

Given any non-constant inner function $\theta$, we have $K_{\theta}\neq \{0\}$. An approach to this result, which gives more information on the structure of model spaces, uses the following factorisation result.

%Theorem 3.2
\begin{thm}\label{lem:16}
Given any non-constant inner function $\theta$,  we may choose $a \in \RR$ and inner functions $\theta_1$, $\theta_2$ where  $\theta_1$ is non-constant, analytic in a neighbourhood of $a$ and $\theta_1(a)=1$,
such that $\theta=\theta_1 \theta_2$.
\end{thm}

\beginpf If $\theta$ has an elementary Blaschke factor $b$, then the result is clear, taking $a=0$ and $\theta_1=b/b(0)$. So we may
assume that $\theta$ is a singular inner function.

If the measure $\mu$ determining $\theta$ is an atom concentrated at
$\infty$, then
we may take $a$ to be any finite point, and the result is clear.

Otherwise, let $I$ be any open interval such that $\mu(\RR \setminus I)>0$, and choose $a \in I$.
Define a decomposition of $\mu$ into positive singular measures by setting $\mu=\mu_1+\mu_2$, where $\mu_1(A)=\mu(A \setminus I)$
and $\mu_2(A)=\mu(A \cap I)$. These determine
inner functions $\theta_1$ and $\theta_2$ with the required properties, and by multiplying them
by unimodular constants, if necessary, we may also assume that $\theta_1(a)=1$.
\endpf
%%%%

It is easy to see that, if $\theta_1$ be a non-constant inner function, analytic in a neighbourhood of a point $a \in \RR $, with $\theta_1(a)=1$, and $\Lambda_{\theta_1,a}$ is the function
\beq \label{2.4}
\Lambda_{\theta_1,a}(\xi)=\frac{\theta_1(\xi)-1}{\xi-a}, \quad \xi \in \RR,
\eeq
then $\Lambda_{\theta_1,a} \in K_{\theta_1}$. If, in addition, $\theta_1$ is a singular inner function then
 $\Lambda_{\theta_1^\mu,a} \in K_{\theta_1}$ for all $\mu\in (0,1]$.

%\bigskip
%\beginpf We have $\Lambda_{\theta_1,a} \in H_p^+$ and $\overline{\theta_1}\,\Lambda_{\theta_1,a} \in H_{p}^-$, so that $\Lambda_{\theta_1,a} \in K_{\theta_1}$. Assume now that $\theta_1$ is a singular inner function. Since $\theta_1$ does not vanish and is analytic in $\CC^+\cup D_{a,\varepsilon}$ for some $\varepsilon >0$, we can define an analytic branch of $\theta_1^{\mu}$ in $\CC^+\cup D_{a,\varepsilon}$ with $\theta_1^{\mu}(a)=1$, and we conclude analogously that $\Lambda_{\theta_1^{\mu},a} \in K_{\theta_1}$.
%\endpf

So, if $\theta$ is a Blaschke product, then $\ds \frac{1}{\xi-\overline {z_+}}\in K_\theta$ for every zero $z_+$ of $\theta$. If $\theta$ is a singular inner function, we can write $\theta=\theta_1 \theta_2$ as in Theorem \ref{lem:16} and
\[
\theta_2\Lambda_{\theta_1^\mu,a} \in K_{\theta} \quad \text{for all } \mu \in (0,1].
\]
Otherwise, $\theta=\alpha B\mathcal{S}$ where $\alpha \in \CC$, $|\alpha|=1$, $B$ is a Blaschke product and $\mathcal{S}$ is a singular inner function, and it is easy to see that $K_\theta \supset K_\mathcal{S}$.\\

In any case, we explicitly see that $K_\theta$ is infinite-dimensional unless $\theta$ is a finite Blaschke product. In the latter case, we can write

\beq \label{2.6}
\theta =h_-r^nh_+, \quad \text{with } h_\pm \in \mathcal{G}H_{\infty}^{\pm}, \, n \in \mathbb N
\eeq
and $K_\theta$ is an $n$-dimensional linear space described by
\beq \label{2.7}
K_\theta = h_+\, \spam \left\{\lambda_+^{-1}r^j:\, j=0,1,\hdots, n-1\right\}=h_+K_{r^n}
\eeq
(recall that $\lambda_\pm(\xi)=\xi \pm i$).

Thus, in the case where $\theta$ is a rational inner function, it is clear from \eqref{2.7} that $K_\theta\subset \lambda_+^{-1}H_{\infty}^+\subset H_{\infty}^+$. 
The question whether $K_\theta\subset H_{\infty}^+$ in other cases is fairly delicate and will be dealt with later in this paper.\\

To have a better understanding of infinite-dimensional model spaces $K_\theta$, it will be useful to characterise some dense subsets. 
While $K_\theta$ may not be itself contained in $H_\infty^+$, there are nevertheless dense subsets of $K_\theta$ contained in $\lambda_+^{-1}H_{\infty}^+$. Indeed, for each $w \in \CC^+$, let
\beq \label{2.8}
k_w(\xi)=\frac{i}{2\pi}\frac{1}{\xi-\overline{w}}, \quad \xi \in \RR,
\eeq
and, given an inner function $\theta$, let $k_w^\theta$ be defined for each $w \in \CC^+$ by
\beq \label{2.9}
k_w^\theta(\xi)=\frac{i}{2\pi}\frac{1-\overline{\theta(w)}\theta(\xi)}{\xi-\overline{w}}=P_\theta k_w(\xi).
\eeq
These are the reproducing kernel functions for $K_{\theta}^2$, but they play the same role in $K_{\theta}^p$ for each $p\in (1,\infty)$, namely
\beq \label{2.10}
\ds \int_{-\infty}^{+\infty} f(x)\overline{k_w^\theta(x)}\,dx=f(w)\quad \text{for all } f \in K_{\theta}^p.
\eeq

Let also $f_k^\theta$ be the functions defined, for each $k\in\mathbb{Z}_0^+$, by
\[
f_k^\theta=\frac{r^k}{\lambda_+}-\frac{a_0+a_1\lambda_+ +\hdots +a_k \lambda_+^k}{\lambda_+^{k+1}}\theta,
\]
where $a_j=(\lambda_-^j \bar \theta)^{(j)}_{(-i)}/j!$ , \,$j=0,1,...,k-1$. As in the case of reproducing kernel functions, these are easily recognisable functions of $K_\theta^p$, providing the following density result.
We have $k^\theta_w, f_k^\theta \in K_{\theta}^p \cap \lambda_+^{-1}H^+_\infty$ for all $w \in \CC^+, k\in\mathbb{Z}_0^+$ and $p \in (1,\infty)$, and
\[
K_{\theta}^p = \clos_{H^+_p} \spam \{ k^\theta_w: w \in \CC^+ \}=\clos_{H^+_p} \spam \{ f_k^\theta: k \in \mathbb{Z}_0^+ \}.
\]

%\beginpf
%The linear span of the rational functions $k_w$\,, $w \in \CC^+$, is dense in $H_p^+$, therefore their images under $P_\theta$ constitute a dense set in $K_{\theta}^p$.

%Analogously, since $\left\{\lambda_+^{-1} r^k\, : k \in \mathbb{Z}_0^+\right\}$ is dense in $H_p^+$ and $f_k^\theta=P_\theta(\lambda_+^{-1} r^k)$, we also have
%$K_{\theta}^p=\clos_{H^+_p} \spam \{f_k^\theta: k \in \mathbb{Z}_0^+ \}.$
% \endpf

%Now we study some relations between model spaces associated to different inner functions.

%Definition 3.3
\begin{defn}
For inner functions $\theta_1$ and $\theta_2$, we write $\theta_2 \preceq \theta_1$ if and only if $\theta_2$ divides $\theta_1$,
in the sense that $\theta_1 = \theta_2 \theta_3 $ for some inner function $\theta_3$.

We also write $\theta_2 \prec \theta_1$ if $\theta_1 = \theta_2 \theta_3 $ for some {\em non-constant\/}
inner function $\theta_3$.
\end{defn}

The results in the next theorem may be considered as generally known; see, for instance, \cite{Nikolski}. 

%Theorem 3.4
\begin{thm}\label{thm:prec}
Let  $\theta_1, \theta_2$ and $\theta_3$ be inner functions. We have, for $p \in (1,\infty)$:
\begin{description}
  \item[(i)] $\theta_2 \preceq \theta_1$ if and only if $K_{\theta_2}^p \subset K_{\theta_1}^p$;
  \item[(ii)] $\theta_2 \prec \theta_1$ if and only if $K_{\theta_2}^p \subsetneq K_{\theta_1}^p$;
  \item[(iii)] $\theta_2 \theta_3 \preceq \theta_1$ if and only if $\theta_3 K_{\theta_2}^p \subset K_{\theta_1}^p$;
 \item[(iv)] $\theta_1 \preceq \theta_3 \implies \theta_1 K_{\theta_2} \subset K_{\theta_3\theta_2}$, where the inclusion is strict if $\theta_1$ is not constant.
\end{description}

\end{thm}

An alternative short proof of (i)--(iii) is provided in Section 5 using the characterisation of maximal functions in a model space instead of the $H_p^+$--$H_q^+$ duality.

For any inner functions $\theta_1$, $\theta_2$ we have
\beq \label{2.16}
 K_{\theta_1} \subset K_{\theta_1\theta_2}, \quad \theta_1 K_{\theta_2} \subset K_{\theta_1\theta_2},
 \eeq
 and the
two subspaces at the left-hand side of these inclusions provide a direct sum decomposition
\beq \label{2.17}
K_{\theta_1\theta_2} = K_{\theta_1} \oplus \theta_1 K_{\theta_2}.
\eeq
%Indeed, $K_{\theta_1} \cap \theta_1 K_{\theta_2}=\{0\}$ and, if $\phi_+ \in K_{\theta_1\theta_2}$, then

%\[\phi_+=\theta_1P^-\overline{\theta_1}\phi_++\theta_1P^+\overline{\theta_1}\phi_+\]
%where the first term on the right hand side belongs to $K_{\theta_1}$ and, in the second term, we have $P^+\overline{\theta_1}\phi_+\in K_{\theta_2}$ because $\overline{\theta_2}(P^+\overline{\theta_1}\phi_+)=\overline{\theta_2}\,\overline{\theta_1}\phi_+-\overline{\theta_2}P^-(\overline{\theta_1}\phi_+) \in H_p^-$.\\

For $p=2$, \eqref{2.17} yields an orthogonal decomposition of $K_{\theta_1\theta_2}$. %As a consequence of Theorem \ref{thm:prec} and \eqref{2.17} we see that, for any inner functions $\theta, \theta_1$, if $K_{\theta_1}\subset K_{\theta}$ then
%\beq \label{2.18}
%K_{\theta} = K_{\theta_1} \oplus \theta_1 K_{\theta/\theta_1}
%\eeq
%and thus every model space has a complement with respect to another model space %containing it, which for $p=2$ is its orthogonal complement.\\
We also have the following.

%Theorem 3.5
\begin{thm} Let $\theta, \theta_1$ be inner functions and let $n,m \in \NN$ with $n \geq m$. Then if
\beq \label{2.19}
\theta^n  \preceq r^m\theta_1
\eeq
then
%\beq \label{2.20}
%\theta^{n-m}  \preceq \theta_1
%\eeq
%and
\beq \label{2.21}
K_{\theta_1} = K_{\theta^s} \oplus \theta^s K_{\theta_1/\theta^s}
\eeq
for any $s \in \mathbb{Z}_0^+$, $s \leq n-m$.
\end{thm}

\beginpf
If \eqref{2.19} holds, let $r^m\theta_1=\theta^n\widetilde\theta$ with $\widetilde\theta$ inner. Then
\beq \label{2.22}
\theta_1 = \theta^{n-m}(r^{-m}\widetilde\theta\theta^m)
\eeq
and, since $\theta_1 \in H_\infty^+$, we must have $\theta(i)=0$ or $\widetilde\theta(i)=0$. In any case, $r^{-m}\widetilde\theta\theta^m$ is an inner function and so \eqref{2.22} implies that $\theta^{n-m}  \preceq \theta_1$. Now \eqref{2.21} follows from Theorem \ref{thm:prec} and \eqref{2.17}.
\endpf

Inner functions and model spaces can be related by an equivalence relation as follows.

%definition 3.6
\begin{defn}
If $\theta_1$ and $\theta_2$ are inner functions, we say that $\theta_1 \sim\theta_2$ if and only if there are functions $h_{\pm} \in \mathcal{G}H_\infty^{\pm}$ such that
\begin{equation}\label{3.22}
\theta_1=h_-\theta_2h_+.
\end{equation}
\end{defn}

It is easy to see that we have $\theta_1=h_-\theta_2h_+$ and $\theta_1=\tilde{h}_-\theta_2\tilde{h}_+$ with $h_{\pm} \in \mathcal{G}H_\infty^{\pm}$, $\tilde{h}_{\pm} \in \mathcal{G}H_\infty^{\pm}$, if and only if $\frac{h_-}{\tilde{h}_-}=\frac{\tilde{h}_+}{h_+}=c \in \mathbb{C}\setminus \{0\}$, and we can choose $h_{\pm}$ in \eqref{3.22} such that $\|h_-\|_\infty=\|h_+\|_\infty=1$.

Moreover, if \eqref{3.22} holds for given $\theta_1$, $\theta_2$, then $h_-\overline {h}_+=h_+^{-1}{(\overline {h_-})^{-1}}$; since the left-hand side represents a function in $H_\infty^-$ and the right-hand side represents a function in  $H_\infty^+$, both are constant and we have

\begin{equation}\label{3.24}
\overline {h_+}=h_-^{-1}c,  \quad \overline{h_-}=h_+^{-1}c^{-1}, \quad \text{with } c \in \mathbb{C}\setminus \{0\}.
\end{equation}

%Definition 3.7
\begin{defn}\label{def:3.11}
If $\theta_1$ and $\theta_2$ are inner functions, we say that $K_{\theta_1}\sim K_{\theta_2}$ if and only if
\begin{equation}\label{3.25}
K_{\theta_1}=h_+K_{\theta_2} \quad \text{with } h_+ \in \mathcal{G}H_\infty^+.
\end{equation}
\end{defn}

It is clear that
\begin{equation}\label{3.26}
\theta_1 \sim \theta_2 \Rightarrow K_{\theta_1}\sim K_{\theta_2}
\end{equation}
since, by \eqref{21.7}, \eqref{21.10} and \eqref{3.24}, if \eqref{3.22} holds then
\[
 K_{\theta_1}=\ker T_{\bar{\theta}_1}=\ker T_{h_-^{-1}\bar{\theta}_2h_+^{-1}}=h_+\ker T_{\bar{\theta}_2}=h_+K_{\theta_2}.
\]

If $\theta$ is a finite Blaschke product, then $K_\theta \sim K_{\tilde \theta}$ if and only if $\tilde \theta$ is also a finite Blaschke product of the same degree. However, model spaces associated with infinite Blaschke products may be equivalent, in the sense of Definition \ref {def:3.11}, to model spaces associated to singular inner functions. In particular, for any singular inner function $\theta$ there exists an infinite Blaschke product $B$ such that
\begin{equation}\label{3.281}
K_\theta \sim K_B.
\end{equation}
In fact the function
\begin{equation}\label{3.29}
B=\frac{\theta-a}{1-\bar{a}\theta}
\end{equation}
is a Blaschke product for all $a$ with $|a|<1$ outside a set of measure zero \cite{duren, Frostman}. Thus any inner function $\theta$ can be factorised as
\begin{equation}\label{3.30}
\theta=h_-Bh_+
\end{equation}
where $B$ is a Blaschke product and $h_\pm \in \mathcal{G}H_\infty^\pm$ with
 \begin{equation}\label{3.31}
h_-=1+a\bar{B}, \quad h_+=\frac{1}{1+\bar{a}B}\,.
\end{equation}
It follows from \eqref{3.30} that $\theta \sim B$ and 
 \begin{equation}\label{3.32}
K_\theta=h_+K_B.
\end{equation}

If $K_{\theta_1}\sim K_{\theta_2}$ then the two model spaces are isomorphic (although not usually isometric in the case $p \ne 2$) and share several properties, namely that they are either both contained in $H_\infty^+$ or they are not.\\
The projections associated with $K_{\theta_1}$ and $K_{\theta_2}$ are related as follows.

%Theorem 3.8
\begin{thm} If $K_{\theta_1}\sim K_{\theta_2}$ and $h_+ \in \mathcal{G}H_\infty^+$ is such that \eqref{3.25} holds, then
\begin{equation}\label{3.27}
\widetilde{P}_{\theta_1}:= h_+P_{\theta_2}h_+^{-1}P^+
\end{equation}
is a projection from $H_p^+$ (or $L_p$) onto $K_{\theta_1}$ such that
\begin{equation}\label{3.28}
\widetilde{P}_{\theta_1}|_{K_{\theta_1}}=P_{\theta_1}.
\end{equation}
\end{thm}

\beginpf $\widetilde{P}_{\theta_1}$ is obviously a projection and, for any $\phi_+ \in H_p^+$, $\widetilde{P}_{\theta_1}\phi_+ \in h_+K_{\theta_2}=K_{\theta_1}$.
Moreover, if $\phi_+ \in K_{\theta_1}$ then $h_+^{-1}\phi_+ \in K_{\theta_2}$, $P_{\theta_2}h_+^{-1}\phi_+=h_+^{-1}\phi_+$, and we have
$\widetilde{P}_{\theta_1}\phi_+=h_+P_{\theta_2}h_+^{-1}\phi_+=\phi_+.$
\endpf

%%%%%%%%%%%%%%%%%%%%%%%%%%%%%%%%%%%%   Section 4     Model spaces contained in Hinfty      %%%%%%%%%%%%%%%%%%%%%%%%%%%%%%%%%%%%%
%%%%%%%%%%%%%%%%%%%%%%%%%%%%%%%%%%%%%%%%%%%%%%%%%%%%%%%%%%%%%%%%%%%%%%%%%%%%%%%%%%%%%%%%%%%%%%%%

\section{Model spaces contained in $H_\infty^+$}
\label{sec:4}

 Let
\beq \label{2.23}
 K_{\theta}^\infty:=H^+_\infty \cap \theta H^-_\infty
\eeq
for an inner function $\theta$. 
Since $\theta \in K_{\theta}^\infty$, we can extend the  inclusion $\theta_1 K_{\theta_2} \subset K_{\theta_1\theta_2}$ in \eqref{2.16} as follows.

\begin{prop}
For any inner functions $\theta_1, \theta_2$ we have
\beq\label{eq:eq9}
K_{\theta_1}^\infty K_{\theta_2} \subset K_{\theta_1 \theta_2}.
\eeq
\end{prop}
\beginpf
Let $f_1^+ \in K_{\theta_1}^\infty$, $f_2^+ \in K_{\theta_2}$. Then $f_1^+f_2^+ \in K_{\theta_1 \theta_2} $ because $f_1^+f_2^+ \in H_{p}^+ $ and
\[\overline{\theta_1}\,\overline{\theta_2}f_1^+f_2^+=(\overline{\theta_1}f_1^+)(\overline{\theta_2}f_2^+) \in H_p^-.\]
\endpf

Using the fact that model spaces are T-kernels and the n.\ $\eta$-invariance of T-kernels for all $\eta \in \overline{K}_\theta= \bar{\theta} K_\theta$ (\cite{CP14}), we also have
\begin{equation}\label{4.3}
K_{\theta_1}K_{\theta_2} \cap H_p^+ \subset K_{\theta_1\theta_2},
\end{equation}
since $K_{\theta_1}K_{\theta_2}=\overline{K}_{\theta_1}(\theta_1K_{\theta_2})$ and $\theta_1K_{\theta_2}\subset K_{\theta_1\theta_2}$ by \eqref{2.16}.

From \eqref{4.3} we have
$
K_{\theta_1}K_{\theta_2}  \subset K_{\theta_1\theta_2}$
if either $K_{\theta_1}$ or $K_{\theta_2}$ is contained in $H_\infty^+$, as happens when $\theta_1$ or $\theta_2$ are finite Blaschke products. The question whether there are infinite-dimensional model spaces satisfying this boundedness condition has different answers depending on  whether the setting is the disk $\mathbb{D}:=\left\{z\in \mathbb{C}:|z|<1\right\}$ or the upper-half plane.

%%%%  Subsection  The case of the disk
\subsection{The case of the disk}

This is the easier case, and the following result holds, which we include for completeness.

%Theorem 4.2
\begin{thm}\label{thm:4.2}
Let $\theta \in H_\infty(\mathbb{D})$ be inner; then, for any $p\in (1,\infty)$, the model space $K_\theta=H_p(\mathbb{D}) \cap \theta \bar{z}\overline{H_p(\mathbb{D})}$ is a subspace of $H_\infty(\mathbb{D})$ if and only if $\theta$ is a rational function.
\end{thm}

\beginpf If $\theta$ is rational then we have $\theta=h_-z^nh_+$ with $h_+,\overline {h_-} \in \mathcal{G}H_\infty(\mathbb{D})$, $\overline{h_\pm}=h_\mp^{-1}$ and $n$ equal to the number of zeros of $\theta$, taking their multiplicity into account. By \eqref{21.7} and \eqref{21.10}, $K_\theta=\ker T_{\bar{\theta}}=h_+\ker T_{\overline{z^n}}$ and it follows that $K_\theta \subset H_\infty(\mathbb{D})$.

Conversely, note that the reproducing kernel functions $k_w^\theta$, with
\[
k_w^\theta(z):=\frac{1-\overline{\theta(w)}\theta(z)}{1-\overline w z},\quad w\in \mathbb{D},
\]
lie in  $K_\theta$, for any $p\in (1,\infty)$.
Indeed, their $H_p(\mathbb{D})$ norm is bounded by a constant times $(1-|w|)^{-1+1/p}$, as can be seen by estimating the norm of
$1/(1-\overline w z)$ directly -- it is enough to consider real positive $w$ and do a direct calculation. This can be
achieved quite simply using an isometry with the $H_p$ space of the half-plane as in \cite[Prop. 2.15]{jrp-hankel}.

However, if $\theta$ is not a finite Blaschke product, then for each $\eps>0$ we can find a point $w\in \mathbb{D}$ with $|w|> 1-\eps$ and
$|\theta(w)|<1/2$. Thus, taking $z=w/|w|$ we have $\|k^\theta_w\|_\infty \ge 1/(2(1-|w|))$, that is
\beq\label{eq:jbdd}
\sup_{f \in K^\theta_p} \frac{ \|f\|_\infty }{ \|f\|_p} = \infty.
\eeq
If every function in $K_\theta$ is bounded then we have a natural embedding $J: K_\theta \to H_\infty(\DD)$. 
%The
%graph of $J$  is closed, since if $f_n \to f$ in $K_\theta$ norm and $J f_n \to g$ in $H_\infty(\DD)$ norm, we clearly have
%$f_n(w) \to f(w)$ and $f_n(w) \to g(w)$ for each $w \in \DD$ (\cite {Nikolski}), meaning that $Jf=g$. 
But the closed graph theorem now implies that $J$ is
a bounded operator, contradicting (\ref{eq:jbdd}).
\endpf

%%% Subsection   The case of the (upper) half-plane

\subsection{The case of the (upper) half-plane}
As in the setting of $H_p$ spaces of the disk, if $\theta$ is a rational inner function then $K^p_\theta \subset H_\infty^+$, for all $p \in (1, \infty)$. Now, however, we may have $K^p_\theta \subset H_\infty^+$ for some classes of irrational inner functions $\theta$, as well as model spaces which are not contained in $H_\infty^+$.

Indeed, Dyakonov \cite{Dy90} (see also \cite{Dy91,Dy00}) gave the following necessary and sufficient
conditions for $K^p_\theta \subset H_\infty^+$ (note that they do not depend on $p$).
\begin{eqnarray}\label{eq:dyk}
1. && \theta' \in H_\infty^+ ;\nonumber \\
2. &&  \inf\{|\theta(z)|: 0 < \imag z < \epsilon\} > 0 \quad \hbox{for some} \quad \epsilon>0.
\end{eqnarray}

In particular, if for $\lambda \in \RR^+$,
$e_\lambda$ denotes the singular inner function
\beq\label{4.6}
e_\lambda(\xi)=e^{i\lambda \xi}, \quad \xi \in \RR,
\eeq
then
for any $p\in (1,\infty)$  the (Paley--Wiener type)
model space $K^p_{e_\lambda}$ consists of entire functions and is contained in $H_\infty^+$.

However, if $\theta$ possesses a sequence of zeroes tending to the real axis, or if $\theta$ has a singular
inner factor other than $e_\lambda$ for some $\lambda>0$, then the model space
$K^p_\theta$ contains unbounded functions. This follows from the well-known fact that, for a singular inner function
determined by a measure $\nu$, the non-tangential boundary
limits are 0 almost everywhere with respect to $\nu$  (see for example \cite[Chap.~1]{duren}).\\

The following result gives an alternative, and occasionally more usable, necessary and sufficient condition for
the inclusion into $H^+_\infty$. 

\begin{thm} \label{k2th}
$K_\theta^p \subset H_\infty^+$ if and only if
\[ 
\sup_{w \in \CC^+} \frac{1-|\theta(w)|^2}{\imag w}<\infty.
\]
\end{thm}

\beginpf
Note that, by Dyakonov's result, it is sufficient to discuss the case $p=2$.
By the closed graph theorem a necessary and sufficient condition for $K^2_\theta$ to embed into
$H^+_\infty$ is that, for all $f \in K^2_\theta$, we have $f\in H_\infty^+$ and
there is a constant $C>0$ such that
\beq\label{eq:4.6a}
\|f\|_\infty \le C \|f\|_2
\eeq
for all $f \in K^2_\theta$.

Since 
 for $f\in K_\theta^2$ we have $\sup_{w \in \CC^+}|f(w)|=\sup_{w \in \CC^+}| \langle f, k_w^\theta\rangle|$, condition \eqref {eq:4.6a} is equivalent to
the condition
that the $L_2$ norms of the $k_w^\theta$ are uniformly bounded, independently of $w$. 
 For $p=2$ we have
\[
\|k_w^\theta\|_2^2 = \langle k_w^\theta, k_w^\theta \rangle = |k_w^\theta(w)| 
\]
and the result follows from  \eqref {2.9}.

\endpf

The following refinement of \eqref{eq:dyk}  is an immediate consequence of Theorem \ref{k2th} and
\eqref{eq:dyk} itself.

%Theorem 4.5
\begin{cor} 
We have $K^p_\theta \subset H_\infty^+$ if and only if
\beq \label {eq:4.6c}
\lim_{\eps \rightarrow 0^+} \inf \{|\theta(z)| : 0 < \imag z < \eps\} = 1.
\eeq
\end{cor}

Dyakonov's condition
that $\inf\{|\theta(z)|: 0 < \imag z < \epsilon\} > 0$ for some $\epsilon>0$ has appeared
elsewhere in the literature, being applied to realization theory \cite{jz02}
and finite-time controllability \cite{JPP}. (The context
is the right half-plane but it is easy to transcribe the results for the upper half-plane.)
In particular, for a Blaschke product with zeroes $\lambda_n=x_n+iy_n$, $n \ge 1$,
the condition is shown in \cite{jz02} to be equivalent to the property that $\inf y_n >0$ and
\[
\sup_{x \in \RR} \sum_{n=1}^\infty \frac{y_n}{y_n^2+(x-x_n)^2} < \infty,
\]
which in turn can be expressed as
 a Carleson measure condition on the measure $\mu:=\sum_{n=1}^\infty y_n \delta_{\lambda_n}$, 
tested on reproducing kernels $k_\lambda$ lying on a horizontal line.\\

A more general question, to which we do not know a complete answer except in the case $p=2$, is to ask when a T-kernel
contains only bounded functions.

%%%%%%%%%%%%%%%%%%%%%%%%%%%%%%%%%%%%%%%%   Section 5    Model spaces and minimal kernels     \to Maximal and minimal functions in model spaces %%%%%%%%%%%%%%%%%%%%%%%%%%%%%%%%%%%%%%%

\section{Maximal and minimal functions in model spaces}
\label{sec:5}

It was shown in \cite{CP14} that for every $\phi_+ \in H_p^+\setminus \{0\}$ there exists a T-kernel containing $\phi_+$, denoted by $\mathcal{K}_{\min}(\phi_+)$, such that for any $g \in L_\infty$ we have 
\beq\label{5.1}
\phi_+ \in \ker T_g \Rightarrow \mathcal{K}_{\min}(\phi_+) \subset \ker T_g
\eeq

and, if $\phi_+=I_+O_+$ is an inner-outer factorisation of $\phi_+$,
\beq\label{5.2}
 \mathcal{K}_{\min}(\phi_+) = \ker T_{\overline{I}_+\overline{O}_+/O_+}.
\eeq

$\mathcal{K}_{\min}(\phi_+)$ is called the minimal kernel for $\phi_+$. It can be shown moreover that a nontrivial, proper, \nei{S^*} subspace $\mathcal{E}$ of $H_p^+$ $(1<p<\infty)$ is a T-kernel if and only if there exists $\phi_+ \in H_p^+$ such that $\mathcal{E}=\mathcal{K}_{\min}(\phi_+)$, i. e., such that $f_+\in \mathcal{E}$ if and only if $f_+\in H_p^+$ and $\overline I_+\frac {\overline O_+}{O_+}f_+\in H_p^- $, where $\phi_+=I_+O_+$ is an inner-outer factorisation of $\phi_+$ (\cite{CP14}).\\

%Definition 5.1

\begin{defn} If $K=\mathcal{K}_{\min}(\phi_+)$, we say that $\phi_+$ is a \emph{maximal function} for $K$.
\end{defn}

Being T-kernels, model spaces are minimal kernels for some of their elements. Given a model space $K_\theta$, it is thus natural to try to characterise the maximal functions for $K_\theta$.\\

We start by remarking that, writing $\theta=\theta_1\theta_2$ as in Theorem \ref{lem:16} and defining $\Lambda_{\theta_1, a}$ as in \eqref{2.4}, we have (for any $p \in \, (1, \infty)$)
\beq\label{5.3}
 \theta_2\Lambda_{\theta_1, a} \in K_{\theta_1\theta_2}=K_\theta.
\eeq
Since $\Lambda_{\theta_1, a}$ is outer,  it follows from \eqref{5.2} that
\beq\label{5.4}
\mathcal{K}_{\min}( \theta_2\Lambda_{\theta_1, a})=\ker T_{\overline{\theta_1 \theta_2}}=K_\theta  .
\eeq

Depending on the inner function $\theta$ associated with the model space, other maximal functions can be defined for $K_\theta$, which may also be useful. The following theorems describe, in different ways, the maximal functions of a given model space $K_\theta$.

%Theor 5.2
\begin{thm}\label{thm:5.2} $K_\theta=\mathcal{K}_{\min}(\phi_+)$ if and only if $\phi_+ \in H_p^+$ and $\phi_+=\theta \phi_-$ with $\phi_-$ outer in $H_p^-$.
\end{thm}

\beginpf If $K_\theta=\mathcal{K}_{\min}(\phi_+)$, then $\phi_+ \in K_\theta$, so that $\phi_+ \in H_p^+$ and 
$\overline\theta\phi_+=\phi_-$ with $\phi_- \in H_p^-$. If $\phi_- $ is not outer in $H_p^-$, then $\phi_-=I_-O_-$ where $I_-$ is a 
non-constant inner function in $H_\infty^-$ and $O_-$ is outer in $H_p^-$. Thus 
$\phi_+ \in \ker T_{\overline{I_-} \overline{\theta}}\varsubsetneq \ker T_{\bar{\theta}}=K_\theta$, which contradicts the assumption.

Conversely, if $\phi_+ \in H_p^+$ and $\phi_+=\theta\phi_-$ with $\phi_-$ outer in $ H_p^-$, then $\phi_+ \in K_\theta$. Moreover, for any $g \in L_\infty$, if $\phi_+ \in \ker T_g$ then $g\theta\phi_- = \eta_- \in H_p^-$, so that
\[g=\bar{\theta}\frac{\eta_-}{\phi_-} \]
where $\phi_-$ is outer in $H_p^-$. Thus, for any $\psi_+ \in H_p^+$ such that $\bar{\theta} \psi_+= \psi_- \in H_p^-$, i.e., for any $\psi_+ \in K_\theta$, we have
\[g\psi_+=\bar{\theta}\frac{\eta_-}{\phi_-}\psi_+=\frac{\eta_-\psi_-}{\phi_-} \in H_p^-\]
because the right-hand side represents a function which is in $L_p$ and in the Smirnov class $\overline{\mathcal{N}_+}$. It follows that
$\psi_+ \in \ker T_g$. Thus $K_\theta \subset \ker T_g$ and we have $K_\theta=\mathcal{K}_{\min}(\phi_+)$.
\endpf

% Remark 5.3
\begin{rem}\label {rem:5.3}

The result of Theorem \ref {thm:5.2} provides an alternative proof to some properties in Theorem \ref{thm:prec} that were
proved using the $L_p-L_q$ duality. Consider, for instance, Theorem \ref{thm:prec} (i) and assume that $K_{\theta_2}^p\subset K_{\theta_1}^p$.
Let $\phi_{\theta_2}^+$ be a maximal function for $K_{\theta_2}^p$, so that by Theorem \ref{thm:5.2} we have $\phi_{\theta_2}^+=\theta_2O_{2-}$
where $O_{2-}$ is outer in $H_p^-$ . Since $K_{\theta_2}^p\subset K_{\theta_1}^p$, then
$\theta_2O_{2-}=\theta_1\psi_-$ with $\psi_-\in H_p^-$; if $\psi_-=I_-O_-$  is an inner-outer factorisation (in $H_p^-$)
then it follows that $\bar{\theta}_2I_-O_-=\bar{\theta}_1O_{2-}$ and, by the uniqueness of inner-outer factorisations,
we conclude that $\bar{\theta}_2I_-=\lambda\bar{\theta}_1$ ($\lambda \in \mathbb{C}$), whence $\theta_2 \preceq \theta_1$.
The same reasoning can be applied to prove (iii) in Theorem \ref{thm:prec}.\\
\end{rem}

As a consequence of Theorem \ref{thm:5.2} we also have:\\

%Theorem 5.4
\begin{thm}\label{thm:5.6} If $\mathcal{K}_{\min}(\phi_+)$ is a model space $K_{\theta_1}$, then $\mathcal{K}_{\min}(\theta \phi_+)$ is also a model space and we have
\beq \label{5.10}
\mathcal{K}_{\min}(\theta\phi_+)=K_\theta \oplus \theta\mathcal{K}_{\min}(\phi_+)=K_{\theta \theta_1} .
\eeq
\end{thm}

\beginpf If $\mathcal{K}_{\min}(\phi_+)=K_{\theta_1}$, where $\theta_1$ is an inner function, then by Theorem \ref{thm:5.2} we have $\phi_+=\theta_1\phi_-$ with $\phi_-$ outer in $H_p^-$. Therefore $\theta\phi_+=\theta\theta_1\phi_-$ and, using Theorem \ref{thm:5.2} again, $\mathcal{K}_{\min}(\theta\phi_+)=K_{\theta\theta_1}$. Since $K_{\theta\theta_1}=K_\theta \oplus \theta K_{\theta_1}$ by \eqref{2.17}, we conclude that \eqref{5.10} holds.
\endpf

We have the following relation for maximal functions in model spaces that are equivalent in the sense of Definition \ref {def:3.11}.

%Theor 5.5

\begin{thm}\label{thm:5.3} Let $\theta_1, \theta_2$ be inner functions and let $K_{\theta_1} \sim K_{\theta_2}$. If \eqref {3.25} holds, then $\phi_+$ is a maximal function for $K_{\theta_1}$ if and only if $\phi_+=h_+\psi_+$, where $\psi_+$ is a maximal function for $K_{\theta_2}$.
\end{thm}

\beginpf Let $\psi_+$ be a maximal function for $K_{\theta_2}$ and let $\psi_+=I_+O_+$ be its inner-outer factorisation. Thus
\[K_{\theta_2}=\mathcal{K}_{\min}(\psi_+)=\ker T_{\overline{I}_+\overline{O}_+/O_+}\]
by \eqref{5.2}. On the other hand, if $\phi_+=h_+\psi_+$ then

\[\mathcal{K}_{\min}(\phi_+) =\mathcal{K}_{\min}(h_+\psi_+) = \ker T_{\overline{I}_+\frac{\bar{h}_+\overline{O}_+}{h_+O_+}} = h_+\ker T_{\overline{I}_+\frac{\overline{O}_+}{O_+}}= h_+K_{\theta_2}\]

by \eqref{21.7} and \eqref{21.10}. Now it follows from \eqref{3.25} that $\mathcal{K}_{\min}(\phi_+) =K_{\theta_1}$. Conversely, if $\phi_+$ is a maximal function for $K_{\theta_1}$ then, from the first part of the proof,
\[\mathcal{K}_{\min}(h_+^{-1}\phi_+)=h_+^{-1}K_{\theta_1}=K_{\theta_2}  \]
and thus $h_+^{-1}\phi_+$ is a maximal function for $K_{\theta_2}$.
\endpf

If $B$ is a Blaschke product vanishing at $z_0^+ \in \mathbb{C}^+$, we have from \eqref{5.2}
\begin{equation}\label{5.6}
K_B=\mathcal{K}_{\min} \left( \frac{B}{\xi-z_0^+} \right).
\end{equation}
Thus it follows from Theorem \ref{thm:5.3} and \eqref {3.26} that if $\theta$ is any non-constant inner function which can be factorised as in \eqref{3.30}, a maximal function for $K_\theta$ will be
\beq\label{5.7}
\phi_+^\theta=h_+\phi_+^B, \quad \text{with } \quad \phi_+^B=\frac{B}{\lambda_{z_0^+}}\ ,
\eeq
where
\beq\label{5.7A}
\lambda_{z_0^+}(\xi):=\xi-z_0^+
\eeq
and we assume that $B(z_0^+)=0$.\\

Note that $\phi_+^\theta$ and $\phi_+^B$ in \eqref{5.7}, as well as the maximal functions in \eqref{5.4}, do not depend on $p$ and belong to $\lambda_+^{-1}H_\infty^+$ (whether or not $K_\theta^p$ is contained in $H_\infty^+$).\\

We can also see that, given any inner function $\theta_1$, from \eqref{3.30} and \eqref{5.7} we have
\beq \label{5.8}
\theta_1=h_-\lambda_{z_0^+}\phi_+^{\theta_1}
\eeq
and that the decomposition $K_{\theta\theta_1}=K_{\theta_1} \oplus \theta_1 K_{\theta}$ (where $\theta$ is an inner function) can also be written in terms of a maximal function for $K_{\theta_1}$ as
\beq \label{5.9}
K_{\theta\theta_1}=K_{\theta_1} \oplus h_-\lambda_{z_0^+}\phi_+^{\theta_1} K_{\theta}\, .
\eeq

Another property relating model spaces with minimal kernels is the following.

%Theorem 5.6
\begin{thm}\label{thm:40}
Let $\phi_1^+,\phi_2^+,\ldots,\phi_n^+$ be such that $\Kmin{\phi_j^+}=K_{\theta_j}$
for each $j=1,2,\ldots,n$, where $\theta_j$ is an inner function.
Then there is a minimal kernel $K$ containing $\{\phi_j^+: j=1,2,\ldots,n\}$, and for
$\theta=\LCM(\theta_1,\theta_2,\ldots,\theta_n)$ we have
\[
K=K_\theta = \clos_{H^+_p}(K_{\theta_1}+\cdots+K_{\theta_n}) = K_{\theta_j} \oplus \theta_j K_{\theta \overline{\theta_j}}
\]
for each $j$.
\end{thm}

\beginpf

$\clos_{H^+_p}(K_{\theta_1}+\cdots+K_{\theta_n})$ is a closed subspace of $H_p^+$, invariant for $S^*=T_{r^{-1}}$, so it is a model space $K_{\widetilde\th}$.
Now $K_{\widetilde\th}$ is a T-kernel, and $K_{\widetilde\th} \supset \{ \phi_1^+,\phi_2^+,\ldots,\phi_n^+ \}$. Since every T-kernel containing
$\{ \phi_1^+,\phi_2^+,\ldots,\phi_n^+ \}$ must be closed, and contain each $K_{\th_j}$, it also contains $K_{\widetilde\th}$, so that the latter is the minimal kernel $K$.

Since $K_{\widetilde\th} \supset K_{\th_j}$, we have $\theta_j \preceq \widetilde\th$, for every $j$, by Theorem \ref{thm:prec} and, since $\th=\LCM(\th_1,\ldots,\th_n)$, we have $\theta \preceq \widetilde\th$. On the other hand,
$K_{\widetilde\th} \subset K_\th$, since $K_{\widetilde\th} \subset H_p^+$ and $\ol\th K_{\widetilde\th} \subset H_p^-$; therefore,
$\widetilde\th \preceq \th$. It follows that $\widetilde\th=\th$.
\endpf

As a motivation for the next definition, we remark now that if $\phi_+=I_+O_+$ is the inner-outer factorisation of a maximal function for $\ker T_g$, so that $\ker T_g=\ker T_{\overline{I}_+\overline{O}_+/O_+}$, it may happen that
\beq \label{5.9.1}
\overline{O}_+/O_+=\overline{I}_{1+}\overline{O}_{1+}/O_{1+},
\eeq
where $I_{1+}$ is a non-constant inner function and $O_{1+}$ is an outer function in $H_p^+$ (take for instance $O_+(\xi)=\frac{1}{(\xi+i)^2}$). In that case, we have
\beq \label{5.9.2}
\ker T_g=\ker T_{\overline{I}_+\frac{\overline{O}_+}{O_+}}=\ker T_{\overline{I_+I_{1+}}\frac{\overline{O}_{1+}}{O_{1+}}},
\eeq
where $I_+ \prec I_+ I_{1+}$. This cannot happen, however, when $\ker T_{\overline{O}_+/O_+}=\spam \{O_+\}$, which is equivalent to saying that $O_{+}^2$ is rigid in $H^+_{p/2}$ (\cite{CP14}). In fact, \eqref{5.9.2} would imply  that
$I_+I_{1+}O_{1+}\in  \ker T_{\overline{I}_+\frac{\overline{O}_+}{O_+}}$ and thus $I_{1+}O_{1+}\in \ker T_{\overline O_+/O_+}=\spam \{O_+\}$, which is impossible for non-constant $\theta_1$.\\

%definition 5.7

\begin{defn}\label{defn:5.4}
If $g \in L_\infty$, we say that $O_+$ is a \emph{minimal function} for $\ker T_g$ if and only if  for some inner function $I_+$ we have $\ker T_g=\mathcal{K}_{\min}(I_+O_+)$ and $\mathcal{K}_{\min}(O_+)=\spam \{O_+\}$ .\\
\end{defn}

In $H_2^+$, every non-trivial T-kernel has a minimal function (\cite{sarason88},\cite{sarason}). The following theorem shows that this property also holds for model spaces in $H_p^+$; whether the same is true in general for T-kernels in $H_p^+$ is an open question, to the authors' knowledge.

%Theor 5.8
\begin{thm}\label{thm:5.5}
For any $p\in (1,\infty)$ and any inner function $\theta$, there exists a minimal function $O_+$ in $K_\theta$.
\end{thm}
\beginpf
With the notation of \eqref{3.30} and \eqref{5.7}, it is enough to consider $O_+=\frac{h_+}{\lambda_{\overline{z}_0^+}}$ and $I_+=B\frac{\lambda_{\overline{z}_0^+}}{\lambda_{z_0^+}}$.
\endpf

%%%     Section 6   On the relations between $\ker T_g$ and $\ker T_{\theta g}$   %%%%%%%%%%%%%%%%%%%% 

\section{On the relations between $\ker T_g$ and $\ker T_{\theta g}$}
\label{sec:6}

If $\theta$ is a non-constant inner function, $g\in L_\infty$ and $\ker T_g \neq \left \{0\right\}$, we have $\ker T_{\theta g} \varsubsetneq \ker T_g$. We may then ask how 
much ``smaller" $\ker T_{\theta g}$ is, with respect to $\ker T_g$, and in particular when is it non-trivial.

%Definition 6.1
\begin{defn} Let $g \in L_\infty$ and $\theta$ be an inner function. If $\ker T_g\neq \{0\}$  and $\ker T_{\theta g}=\{0\}$, we say that \emph{$\theta$ annihilates $\ker T_g$.}
\end{defn}

It is clear that a necessary and sufficient condition for $\ker T_g$ not to be annihilated by $\theta$ is that there exists $\phi_+$ such that
\beq\label{6.3}
 \theta \phi_+ \in \ker T_g, \ \ \ \phi_+ \in H_p^+\setminus\{0\},
\eeq
and in this case $\phi_+ \in \ker T_{\theta g}$.\\

If $\theta$ is a finite Blaschke product we have the following result from \cite{BCD}, taking into account that in this case $\theta\sim r^k$, where $k$ is the number of zeroes of $\theta$.

%Theorem 6.2
\begin{thm}\label{thm:6.3}
If $g \in L_\infty$ and $\theta$ is a finite Blaschke product, then
\beq\label{6.4}
\dim \ker T_g < \infty \Leftrightarrow \dim \ker T_{\theta g} < \infty.
\eeq
We have $\dim \ker T_g < \infty$ if and only if there exists $k_0   \in \ZZ$ such that $\ker T_{r^{k_0}g}=\{0\}$ and, in this case, $\dim \ker T_g \leq \max \{0,k_0\}$. Moreover, if $\dim \ker T_g < \infty$, we have
\beq\label{6.5}
\dim \ker T_{\theta g} =\max \{0, \dim \ker T_g - k \}
\eeq
 where $k$ is the number of zeroes of $\theta$ counting their multiplicity.
\end{thm}

 Thus, in particular, if $\dim \ker T_g=d< \infty$ and $\theta$ is a finite Blaschke product such that $\dim K_\theta \leq d$, then
 \[
 \dim \ker T_{\theta g} = \dim \ker T_g - \dim K_\theta\,.
 \]
   If $\theta$ is not a finite Blaschke product and $\dim\ker T_g < \infty$, then $\ker T_{\theta g}=\{0\}$, since $\theta \phi_+ \in \ker T_g$ implies that $\theta_1 \phi_+ \in \ker T_g$ for all inner function $\theta_1$ such that $\theta_1 \prec \theta$. On the contrary, if $\ker T_g$ is infinite-dimensional then
$\ker T_{\theta g}$ may or may not be finite-dimensional, and in particular it may be $\{0\}$. It is clear that $\theta$ annihilates $\ker T_g$ if $\bar{g}\in H_\infty^+$ is an inner function and $\theta \succ \bar{g}$, but that may also happen when no such relation holds between $\theta$ and $\bar{g}$, as in the example that follows.

%Example 6.3
\begin{ex}
Let $g(\xi)=e^{i /\xi}$, $\theta(\xi)=e^{i\xi}$. For $p=2$, we have

\beq\label{6.1}
f_+ \in \ker T_{\theta g} \Leftrightarrow f_+ \in H_2^+\, ,\, e^{i\xi}e^{i /\xi}f_+=f_- \in H_2^-.
\eeq

Using the isometry from $H_2^+$ onto $H_2^-$ defined by $f\mapsto \tilde{f}$ with $\tilde{f}(\xi)=\frac{1}{\xi}f(\frac{1}{\xi})$, we obtain from \eqref{6.1}:
\beq\label{6.2}
 e^{i\xi}e^{i /\xi}f_+=f_-  \Leftrightarrow e^{-i\xi}e^{-i /\xi}\tilde{f}_-=\tilde{f}_+
\eeq
$(\tilde{f}_\pm\in H_2^\mp) $. Since, by Coburn's Lemma, we have $\ker T_{\theta g}=0$ or $\ker T_{\overline{\theta g}}=0$, it follows from \eqref{6.2} that $f_+=0$. 
Therefore, in this case, $\ker T_g$ is infinite-dimensional and $\ker T_{\theta g}=\{0\}$.
\end{ex}

Condition \eqref{6.3} implies a certain ``lower bound" for T-kernels not to be annihilated by an inner function $\theta$. We have the following.

%%Theorem 6.4
\begin{thm}\label{theor:6.41}
Let $g \in L_\infty$ and $\theta\in H_\infty^+$ be an inner function. 
Suppose that $\ker T_{\theta g}\neq \{0\}$, and let $\phi_+$ be a maximal function for $\ker T_{\theta g}$. Then, for any $z_0\in \CC^+$ and any $h_- \in \mathcal{G}H_\infty^-$,
\beq\label{6.6}
\ker T_{g} \supset (h_-\lambda_{z_0}\phi_+K_\theta \cap H_p^+)\oplus \ker T_{\theta g}
\eeq
where $\lambda_{z_0}(\xi)=\xi-z_0$.
\end{thm}

\beginpf
We have $K_\theta =\theta\, \overline{K_\theta}$ with $\overline{K_\theta}\subset \widetilde{\mathcal{N}}_p$ and we also have $h_-, \lambda_{z_0} \in \widetilde{\mathcal{N}}_p$. Thus if $\phi_+ \in \ker T_{\theta g}$, which is equivalent to $\theta\phi_+ \in \ker T_{g}$, it follows that $h_-\lambda_{z_0}\bar{k}_+\theta\phi_+ \in \ker T_g$ for all $k_+ \in K_\theta$ such that the left-hand side of this relation represents a function in $H_p^+$. Thus $(h_-\lambda_{z_0}\phi_+K_\theta \cap H_p^+) \subset \ker T_g$.\\
Clearly, we also have $\ker T_{\theta g}\subset \ker T_g$. Moreover, as we show next,
\beq\label{6.7}
h_-\lambda_{z_0}\phi_+K_\theta \cap \ker T_{\theta g} = \{0\}.
\eeq
To prove this, we start by remarking that $\ker T_{\theta g}=\ker T_{h_-^{-1}\theta g}$. Now assume that $\mathcal{K}_{\min}(\phi_+)=\ker T_{\theta g}$ and $\phi_+=I_+O_+$ is an inner-outer factorisation; let moreover $\psi_+=h_-\lambda_{z_0}\phi_+k_+$, with $k_+ \in K_\theta$, be a function in $H_p^+$. Then
\[\psi_+ \in \ker T_{\theta g}\Leftrightarrow \psi_+ \in \ker T_{\overline{I_+} \overline{O_+}/O_+} = \ker T_{h_-^{-1}\overline{I_+} \overline{O_+}/O_+}\]
\[\Leftrightarrow \lambda_{z_0}k_+\overline{O_+}=\psi_- \in H_p^-.\]
Therefore we have $k_+=\frac{\psi_-}{\overline{O_+}\lambda_{z_0}} \in \overline{\mathcal{N}}_+\cap L_p =H_p^-$ and, since $k_+ \in H_p^+$, it follows that $k_+=0$. Thus
\[ (h_-\lambda_{z_0}\phi_+K_\theta \cap H_p^+)\cap \ker T_{\theta g}=h_-\lambda_{z_0}\phi_+K_\theta\cap \ker T_{\theta g}=\{0\}.\]
\endpf

%Remark 6.5
\begin{rem}
Let $h_-=1$, $z_0^+=i$ (so that $\lambda_{z_0^+}=\lambda_-$) and let $f=\lambda_-\phi_+$, $\mathcal{K}=\spam \{P_\theta(\lambda_+^{-1}r^k): k \in \ZZ_0^+\}$. The previous result implies that whenever $\ker T_{\theta g}\neq \{0\}$ we must have
\beq\label{6.9}
\ker T_g \supset f \mathcal{K} \oplus \ker T_{\theta g},
\eeq
where $f\neq 0$ and $\mathcal{K}$ is dense in $K_\theta$.
\end{rem}

Moreover, with the same assumptions as in Theorem \ref {theor:6.41}:

%Corollary 6.6
\begin{cor} \label{cor:6.6}
If $h_-\lambda_{z_0}\phi_+K_\theta \subset H_p^+$ then, for  $f=h_-\lambda_{z_0}\phi_+$, we have $\ker T_g \supset f K_\theta \oplus \ker T_{\theta g}$ .
\end{cor}

In particular, if $\theta g=\bar \theta_1$, then $\ker T_{\theta g}$ is a model space $K_{\theta_1}$, and $\ker T_g =K_{\theta \theta_1}$. Choosing for $K_{\theta_1}$ a maximal function $\phi_+^{\theta_1}$ such that $\theta_1=h_-\lambda_{z_0} \phi_+^{\theta_1}$ as in \eqref{5.8}, we see from \eqref{5.9} that the inclusion in Corollary \ref {cor:6.6} becomes an equality in this case.

Another case in which the inclusions of Theorem \ref{theor:6.41} and Corollary \ref{cor:6.6} can also be replaced by equalities is the one that we study below.
\\

We start by remarking that, in the case of an infinite-dimensional $\ker T_g$, it follows from Theorem \ref{thm:6.3} that, if $\theta$ is a finite Blaschke product, then $\ker T_{\theta g}$ is an infinite-dimensional proper subspace of $\ker T_g$. Thus it is not possible to relate their dimensions as in Theorem \ref{thm:6.3} for finite-dimensional T-kernels. We can, however, present an alternative relation which not only generalises Theorem \ref{thm:6.3} but moreover sheds new light on the meaning of \eqref{6.5} when $k < \dim \ker T_g < \infty$.

\vspace{0,2cm}

Let $r_z(\xi):=\frac{\xi-z}{\xi-\bar{z}}$ and let
\[B=B_1\cdot B_2  \cdots B_n\]
with $B_j=r_{z_j}^{k_j}$, $j=1,2,\hdots n$, and $k_j \in \NN$, $z_j \in \CC^+$ for each $j=1,2,\hdots n$.

Let moreover
\[k=\sum_{j=1}^{n}k_j.\]

With this notation, we have the following.

%Thorem 6.7
\begin{thm}\label{thm:6.7}
Let $g \in L_\infty$. If $\dim \ker T_g\leq k$, then $\ker T_{Bg}=\{0\}$; if $\dim \ker T_g>k$, then
\beq\label{6.10}
\ker T_g=\ker T_{Bg}\oplus \lambda_{z_1}\phi_+K_{B}
\eeq
where
\beq\label{6.10A}
\lambda_{z_1}(\xi):=\xi-z_1
\eeq
and $\phi_+$ is a maximal function for $\ker T_{Bg}$, i.e.,
\beq\label{6.10B}
\mathcal{K}_{\min}(\phi_+)=\ker  T_{Bg}.
\eeq
\end{thm}

\beginpf
If $\dim \ker T_g >k$, then $\ker T_{Bg}\neq \{0\}$ by Theorem \ref{thm:6.3};  let $\phi_+$ be a maximal function for $\ker T_{Bg}$. Since, for any inner function $\theta \in H_\infty^+$, $\mathcal{K}_{\min}(\phi_+)=\ker T_{\theta g}$ implies that $\tilde{\theta} \phi_+ \not\in \ker T_{\theta g}$ whenever $\tilde{\theta}$ is a non-constant inner function, we have that
\beq\label{6.11}
\widetilde{B}\phi_+ \in \ker T_g \setminus \ker T_{B g} \quad \text{if } \widetilde{B}\preceq B, \ \widetilde{B} \not\in \CC.
\eeq
Let us define, for $g \in L_\infty$,
\beq\label{6.12}
 (\ker T_{g})_-:=g \ker T_g\subset H_p^-.
\eeq
It is easy to see that $(\ker T_{g})_-$ is nearly $\alpha_+$-invariant for all $\alpha_+ \in H_\infty^+$, in the sense that
\beq\label{6.13}
\alpha_+(\ker T_{g})_- \cap H_p^- \subset (\ker T_{g})_-.
\eeq
Let moreover
\beq\label{6.14}
\phi_-=gB\phi_+.
\eeq
It is clear that $\phi_-$ cannot have a non-constant inner factor (in $H_\infty^-$), i.e., $\phi_-$ is an outer function in $H_p^-$; otherwise there would be some non-constant inner function $\theta \in H_\infty^+$ such that $\phi_-=\bar{\theta} \tilde{\phi}_-$ with $\tilde{\phi}_- \in H_p^-$, and it would follow from \eqref{6.14} that $\phi_+ \in \ker T_{\theta B g}\varsubsetneq \ker T_{B g}$, contradicting \eqref{6.10B}. Therefore,

\beq\label{6.15}
\phi_-(\bar{z}_j)\neq 0, \quad \text{for all } j=1,2, \hdots, n.
\eeq

From \eqref{6.14}, \eqref{6.15} and \eqref{6.11} we also see that not only
\beq\label{6.16}
\phi_-\in (\ker T_g)_- \setminus \overline{B} H_p^-
\eeq
but also
\beq\label{6.17}
\bar{\beta}\phi_-\in (\ker T_g)_- \setminus \overline{B} H_p^- \quad \text{if } \beta \prec B,
\eeq
where $\beta$ is an inner function.

Let now $\psi_-$ be any element of $(\ker T_g)_-$. We have
\beq\label{6.18}
\psi_- - \frac{\psi_-(\bar{z}_1)}{\phi_-(\bar{z}_1)}\phi_-=r^{-1}_{z_1}\tilde{\psi}_{1-} \in (\ker T_g)_-
\eeq
where, by \eqref{6.13}, $\tilde{\psi}_{1-} \in (\ker T_g)_-$. Repeating the same reasoning $k_1$ times, we get (for some constants $a_0, a_1, \hdots , a_{k_1-1}$),

\begin{eqnarray}
\psi_- &=& (a_0+ a_1r^{-1}_{z_1}+ \hdots + a_{k_1-1}r^{-(k_1-1)}_{z_1})\phi_- + \overline{B}_1 \psi_{1-} \nonumber \\
 &=& p_{z_1}^- \phi_- + \overline{B}_1 \psi_{1-} \label{6.19}
\end{eqnarray}

where $p_{z_1}^- \phi_- \in (\ker T_g)_- \setminus  \overline{B}H_p^-$ by \eqref{6.17}, $\overline{B}_1 \psi_{1-} \in (\ker T_g)_-$ and $\psi_{1-} \in (\ker T_g)_-$ by \eqref{6.13}.

Analogously, for some constants $b_0, b_1, \hdots, b_{k_2-1}$, we have

\begin{eqnarray}
\psi_- &=& (b_0+ b_1r^{-1}_{z_2}+ \hdots + b_{k_2-1}r^{-(k_2-1)}_{z_2})\phi_- + \overline{B}_2 \psi_{2-} \nonumber \\
 &=& p_{z_2}^- \phi_- + \overline{B}_2 \psi_{2-} \label{6.20}
\end{eqnarray}
and substituting in \eqref{6.19} we obtain
\beq\label{6.21}
\psi_- = (p_{z_1}^-  + \overline{B}_1 p_{z_2}^-)\phi_{-} + \overline{B}_1 \overline{B}_2 \psi_{2-}
\eeq
with
\beq\label{6.22}
(p_{z_1}^-  + \overline{B}_1 p_{z_2}^-)\phi_{-} \, \in \, (\ker T_g)_- \setminus \overline{B}H_p^-
\eeq
\beq\label{6.23}
 \overline{B}_1 \overline{B}_2 \psi_{2-} \, \in \, (\ker T_g)_-.
\eeq
Assuming, for simplicity, that $n=2$, \eqref{6.23} is equivalent to 
\beq\label{6.24}
  \overline{B} \psi_{2-} \, \in \, (\ker T_g)_-\cap \overline{B}H_p^-.
\eeq
Since
\[(p_{z_1}^-  + \overline{B}_1 p_{z_2}^-)\phi_{-} \, \in \, \lambda_{z_1} \phi_- \overline{K}_B= \lambda_{z_1} \phi_- \overline{B}K_B\]
and
\[(\ker T_g)_-\cap \overline{B}H_p^-=\overline{B}(\ker T_{Bg})_-,\]
it follows from \eqref{6.21}, \eqref{6.22} and \eqref{6.24} that
\[(\ker T_g)_-=\overline{B}(\ker T_{Bg})_-\oplus \lambda_{z_1} \phi_- \overline{B}K_B.\]
Therefore
\[g^{-1}(\ker T_g)_-=B^{-1}g^{-1}(\ker T_{Bg})_- \oplus (B^{-1}g^{-1}\phi_-)\lambda_{z_1} K_B\]
\[\Leftrightarrow \ker T_g = \ker T_{Bg} \oplus \phi_+\lambda_{z_1} K_B.\] \endpf

%% Remark 6.8
\begin{rem}
It is not difficult to see, using the n. $\eta$-invariance of T-kernels for $\eta\in H_\infty^-$, that the decomposition \eqref{6.10} still holds if we replace $\lambda_{z_1}\phi_+K_{B}$ by $h_-\lambda_{z_1}\phi_+K_{B}$, for any $h_-\in H_\infty^-$ such that the latter is contained in $H_p^+$, as happens in \eqref{5.9} for model spaces. %%%%%%%%%%%%%%%%%%%%%%%%%%%%%%%%%%%%%%%%%%%%%%%%%%%%%%%%%%%%%%%%%%%%%%%%%%%%%%%%%%%%%%%%%%%%%%%%%%%%%%%%%%%%%%%%
For $p=2$, we may ask whether, by choosing appropriate functions $\phi_+$ and $h_-$ as in \eqref{5.9}, we can make the direct sum in \eqref{6.10} orthogonal. 
\end{rem}

Theorem \ref{thm:5.6} implies that if $\phi_+$ is a maximal function for a model space $K_{\theta_1}=\ker T_{\overline\theta_1}$, then $\theta \phi_+$ is a maximal function for the model space $K_{\theta\theta_1}=\ker T_{\overline{\theta\theta_1}}$ (where $\theta$ is any inner function). As a consequence of Theorem \ref {thm:6.7} we can now generalise this result, when $\theta$ is a finite Blaschke product, to any T-kernel.

%%theorem 6.9

\begin{thm}\label{thm6.8} Let $B$ be a finite Blaschke product and let $g \in L_\infty$. If $\phi_+$ is a maximal function for $\ker T_g$, then $B\phi_+$ is a maximal function for $\ker T_{\overline{B}g}$.
\end{thm}
\beginpf
Assume that $B$ is a (non-constant) finite Blaschke product and let $z_1$ be one of its zeroes. Assume moreover that $\phi_+$ is such that
\[\mathcal{K}_{\min} (\phi_+)=\ker T_g \]
and let $\phi_+=I_+O_+$ be an inner-outer factorisation. Then, by \eqref{5.2},
\[\ker T_g= \ker T_{\overline{I}_+\overline{O}_+/O_+} \quad \text 
{and}\quad \mathcal{K}_{\min} (B\phi_+)=\ker T_{\overline{B}\,\overline{I}_+\overline{O}_+/O_+}. \]
So, from Theorem \ref{thm:6.7},
\begin{eqnarray*}
\mathcal{K}_{\min} (B\phi_+)&=&\ker T_{\overline{I}_+\overline{O}_+/O_+} \oplus \lambda_{z_1}\phi_+ K_B \\
&=& \ker T_g \oplus \lambda_{z_1}\phi_+ K_B = \ker T_{\overline{B}g}.
\end{eqnarray*}
\endpf

%%%%%%%%%%%%%%%%%%%%%%%%%%%%%%%%%%%%%%%%%%%

\section*{Acknowledgments}

This work was partially supported by FCT/Portugal through the projects
PTDC/MAT/121837/2010\,, \, PEst-OE/EEI/LA0009/2013\,, \,and\, PEst-OE/\\MAT/UI0013/2014.

\end{document}